\input amstex
\documentstyle{amsppt}
\magnification=\magstep1
\vcorrection{-0.8cm}
\nologo
\def\ov{\overline}
\def\ovee{\omega^\vee}
\def\Vol{\text{\it Vol}}
\def\des{\text{\it Des\,}}
\def\({\left(}
\def\){\right)}

\def\endemo{\qed\enddemo}
\def\ua{\widehat{\frak{u}}}
\def\w{\wedge}
\def\g{\frak g}
\def\ga{\widehat{\frak g}}
\def\h{\frak h}
\def\ha{\widehat{\frak h}}

\def\ppa{\widehat{\frak p}}
\def\sp{\text{\it Span\,}}

\def\bb{\frak b}
\def\D{\Delta}
\def\l{\lambda}
\def\Dp{\Delta^+}

\def\Da{\widehat\Delta}
\def\Pia{\widehat\Pi}
\def\Dap{\widehat\Delta^+}
\def\Wa{\widehat{W}}
\def\d{\delta}

\def\Ha{\widehat H}

\def\a{\alpha}
\def\b{\beta}

\def\l{\lambda}
\def\d{\delta}
\def\dd{\partial}

\def\th{\theta}

\def\i{{\frak i}}

\def\ganz{\Bbb Z}

\def\real{\Bbb R}

\TagsOnRight

\def\s{\sigma}
\def\l{\lambda}
\def\L{\Lambda}
\def\({\left(} 
\def\){\right)}

\def\g{\frak g}
\def\u{\frak u}

\def\ha{\widehat{\frak h}}
\def\hard{\ha_\real^*}
\def\hauno{\ha^*_1}
\def\hazero{\ha^*_0}
\def\huno{\h^*_1}
\def\hzero{\h^*_0}

\def\C{\Bbb C}
\def\R{\Bbb R}

\def\si{\sigma}
\def\Si{\Sigma}
\def\G{\Gamma}
\def\What{\widehat W}

\def\supp{\text{\it Supp\,}}
\def\ni{\noindent}

\topmatter
\title Abelian subalgebras in $\Bbb Z_2$-graded Lie algebras and affine Weyl groups.
\endtitle
\author Paola Cellini \\ Pierluigi M\"oseneder Frajria\\ Paolo Papi
\endauthor
\address Paola Cellini \vskip 0.pt Dipartimento di Scienze
      \vskip 0.pt Universit\`a  di Chieti-Pescara \vskip 0.pt Viale Pindaro 42
\vskip 0.pt 65127 Pescara --- ITALY \vskip 0.pt e-mail:{\rm \
cellini\@sci.unich.it }
\endaddress
\abstract Let $\g=\g_0\oplus \g_1$ be a simple
$\Bbb Z_2$-graded Lie algebra and let $\frak b_0$ be a fixed Borel subalgebra of 
$\frak g_0$. We  describe and enumerate the abelian $\frak b_0$-stable subalgebras
of $\g_1$.\endabstract
\bigskip
\bigskip
\address Pierluigi M\"oseneder Frajria \vskip 0.pt  Politecnico di Milano\vskip
0.pt Polo regionale di Como\vskip 0.pt Via Valleggio 11
\vskip 0.pt  22100 COMO --- ITALY \vskip 0.pt e-mail:{\rm \
frajria\@mate.polimi.it }
\endaddress
\bigskip
\bigskip
\address Paolo Papi
\vskip 0.pt Dipartimento di Matematica Istituto G. Castelnuovo \vskip 0.pt
Universit\`a di Roma "La Sapienza" \vskip 0.pt Piazzale Aldo Moro 5 \vskip 0.pt
00185 Rome --- ITALY \vskip 0.pt e-mail:{\rm \ papi\@mat.uniroma1.it}
\endaddress
\leftheadtext{Abelian subalgebras in $\Bbb Z_2$-graded Lie algebras}
\rightheadtext{Abelian subalgebras in $\Bbb Z_2$-graded Lie algebras}
\keywords
    abelian ideal, affine Weyl group,  Lie algebra automorphism
\endkeywords
\thanks 
The first and third authors were partially supported by EU grant
HPRN-CT-2001-00272, Algebraic Combinatorics in Europe
\endthanks

\subjclass Primary: 17B20
\endsubjclass

\endtopmatter
\document
\heading\S1 Introduction\endheading\bigskip In this paper we solve the following problem,
which has been posed by D. Panyushev in \cite{12, \S3}. {\sl Suppose that
$\g=\g_0\oplus\g_1$ is a simple
$\Bbb Z_2$-graded Lie algebra and let $\frak b_0$ be a fixed Borel subalgebra of 
$\frak g_0$. Describe and enumerate the abelian $\frak b_0$-stable subalgebras
of
$\g_1$}. We obtain uniform formulas (which will be displayed at the end
of the Introduction) in terms of combinatorial data
associated to the $\Bbb Z_2$-gradation. The interest in this
question lies in  a theorem by Kostant
\cite{7}
 (which has been generalized to the $\Bbb Z_2$ setting by Panyushev \cite{11}) 
relating commutative subalgebras to the maximal eigenvalue of the Casimir
element. More precisely, if $\frak a \subseteq \g_1$ is an abelian $\frak
b_0$-stable subalgebra of dimension $k$,  the corresponding decomposable
$k$-vector in $\Lambda^k\g_1$ obtained by wedging the vectors of a basis of
$\frak a$ is an eigenvector of maximal eigenvalue for the Casimir operator
$\Omega_0$ of $\g_0$. Viceversa, any decompo\-sable element in the ``maximal" eigenspace 
of $\Omega_0$ which is a highest weight vector for the action of $\g_0$ corresponds to an
abelian $\frak b_0$-stable subalgebra.\par Panyushev has solved the previous problem in the very
special case of the little adjoint module (i.e., when $\g_1$ is the irreducible
$\g_0$-module of highest weight $\th_s$, the highest short root of $\g_0$). Panyushev's
strategy consists in identifying, in these cases, the abelian $\frak b_0$-stable subalgebras
of $\g_1$ with the abelian
ideals  consisting only  of long roots of a Borel subalgebra of the Langlands
dual $\g_0^\vee$ of $\g_0$. Then the enumerative result follows
by providing a bijection between these ideals and the  alcoves in the intersection 
of the fundamental chamber of the Weyl group of $\g_0$ with the half-space $(\theta_s,x)<1$.
\par Before describing  our approach to the
general case, let us discuss for a moment Peterson's $2^{rank}$ abelian ideals theorem and its
interpretations, since this result is crucial for our goals. In \cite{8} Kostant attributed to D. Peterson the following
result: the abelian ideals of a Borel subalgebra of a simple Lie algebra $\g$ of rank
$n$ are $2^ n$ in number. Moreover they are parameterized by a special subset
of   the affine Weyl group
$\widehat W$ of $\g$, the one consisting of the so-called minuscule elements. In
\cite{3} a geometric interpretation of the minuscule elements was proposed:
they are exactly the elements which map the fundamental alcove $C_1$ of 
$\widehat W$ into $2C_1$. In particular they are $2^n$ in number, since the ratio
between the volumes of $2C_1$ and $C_1$ is $2^n$. Very recently an alternative
approach to Peterson's result has been proposed in
\cite{9}: identify
$\Lambda^k(\g)$ with
$\Lambda^{(k,k)} \u$, where $\u=x\g[x]$ is the space of $\g$-valued polynomial maps
without constant term and $\Lambda^{(k,k)} \u$ is the subspace of   elements in $\Lambda^k\u$ having $x$-degree $k$. 
Then one can identify $\u$ with the nilradical $\u_F^-$ of the opposite  parabolic subalgebra
$\frak p_F$ corresponding to $\g$ in extended loop algebra of $\g$.
 Now Garland-Lepowsky generalization of Kostant theorem on the cohomology of
$\u_{F}^{-}$ relates the abelian ideals of a Borel subalgebra of $\g$ to the
minuscule elements of $\widehat W$.\par Our approach to the description of the
abelian $\frak b_0$-stable subalgebras of
$\g_1$ is based on a suitable combination of the two ideas described above.
Recall that $\Bbb Z_2$-gradings are in bijection with involutions. Given any
involution
$\sigma$ of $\g$, we give the notion of {\it $\sigma$-minuscule element} in
$\widehat W$, and we prove that  the set $\Cal W_{ab}^\sigma$ of $\sigma$-minuscule
elements is a parameter space for   the abelian
$\frak b_0$-stable subalgebras of
$\g_1$. Then we describe a canonical polytope $D_\sigma$ such that  the
cardinality of $\sigma$-minuscule elements equals
$\frac{Vol(D_\sigma)}{Vol(C_1)}$. Finally we calculate $Vol(D_\sigma)$. 
To describe more in detail the  final outcome we have to fix some 
notation and to recall Kac's classification of involutions of simple Lie
algebras.  
\par 
Let $\g$ be a simple Lie algebra of type $X_N$, $\sigma$ be
an involution of $\g$ and $\g=\g_0\oplus \g_1$ the corresponding
gradation. Let $k$ be the minimal integer such that $\sigma^k$ is
of inner type and consider a realization $(\ha, \Pia, \Pia^\vee)$ of
the affine Kac-Moody algebra $\ga$ of type $X_N^{(k)}$. Set
$\Pia=\{\a_0,\dots,\a_n\}$, and let $a_0, \dots, a_n$ be the labels of
the  affine Dynkin diagram of  $\ga$.
Kac's classification of involutions states that 
$\s$ is completely determined by a $(n+2)$-tuple
$(s_0,\ldots,s_n; k)$, where $n=rk(\g_0)$, $k$ is as  above,  
and $s_0,\ldots,s_n$ are coprime non negative integers satisfying
$k\sum\limits_{i=0}^na_is_i=2$.
We can define a $\ganz$-grading $\ga=\bigoplus\widehat\g_j$ of
$\ga$ (see \S2.1), determined by $\s$,  such that
the degree zero space $\ga_0$ is a
reductive subalgebra of $\ga$ containing $\ha$. We denote by $\Da_0$
the set of roots of $\ga_0$ relative to $\ha$; if $\Dap$ is
the positive system of $\ga$ corresponding to $\Pia$ then we can
choose $\Dap_0=\Dap\cap\Da_0$ as a positive system for $\ga_0$. 
We can view $\g_0$ and $\g_1$ inside $\ga$ as follows. 
Following \cite{6, Ch. 8}, we consider the extended loop algebra 
$ \widehat L(\g,\sigma)=\bigoplus\limits_{j\in\ganz}(\g_{j\,\text{mod}\,
2}\otimes t^{j})\oplus\C c\oplus\C d $ (see \S2.2), with its
natural $\Bbb Z$-grading. Fix a Cartan subalgebra $\h_\s$ in $\g_0$.
There exists an isomorphism of graded algebras $\Psi:
\widehat L(\g,\sigma)\to\ga$ mapping $\h_\s\otimes 1$ into $\ha$.  
Indeed, we can identify the set of
$\h_\s$-roots of $\g_0$ with $\Da_0$, so that $\Dap_0$ defines a
Borel subalgebra $\frak b_0$ of $\g_0$. Moreover, $\Psi$ maps
$\g_1\otimes t^{-1}$ onto $\ga_{-1}$, and  the $\frak b_0$-stable abelian
subalgebras of
$\g_1$ correspond under $\Psi$ to the $\Psi(\frak b_0)$-stable abelian
subalgebras of $\ga_{-1}$.
\par
Now remark that for the $(n+2)$-tuple $(s_0,\ldots,s_n;
k)$ characterizing $\s$ there are only three kind of possibilities.
\roster
\item
$k=1$, and  there exist two indices $p,q$ such that $a_p=a_q=s_p=s_q=1$ 
and
$s_i=0$ for $i\ne p,q$.
\item
$k=1$, and  there exists an   index $p$ such that $s_p=1,\,a_p=2$ and 
$s_i=0$
for $i\ne p$.
\item $k=2$, and  there exists an   index $p$ such that $s_p=1,\,a_p=1$ 
and
$s_i=0$ for $i\ne p$.\endroster 

Denote 
 by $W_\s$ the  Weyl group of $\Da_0$.  Let    $W_f$ be the  Weyl group 
of the root
system generated by
$\Pi_f=\{\a_1,\ldots,\a_n\}$. 
\medskip Consider now case (1),  the hermitian symmetric case. We can assume that
$p=0$, hence we may regard $W_\s$ as a subgroup of $W_f$. Denote by $\ell_\sigma,\,\ell_f$ the  connection indices of
$W_\s,\,W_f$, respectively. We have
$$|\Cal W_{ab}^\sigma|=\frac{|W_f|}{|W_\s|}\left(1+\frac{\ell_\s}{\ell_f}\right).
$$ 
This formula is proved by providing  a  decomposition $D_\s=D_\s'\cup D_\s''$ of the polytope $D_\s$ such that the 
 volumes of $D_\s',\,D_\s''$ are the two summands in the r.h.s of the previous formula. Indeed we prove directly that
$D'_\s$ consists of $\frac{|W_f|\cdot \ell_\si}{|W_\s|\cdot \ell_f}$ alcoves, whereas 
$D_\s''$ is shown to be isometric to a  region whose alcoves are indexed by a set of coset representatives of $W_\s$ in
$W_f$.
\medskip
Now we deal with   cases (2), (3). We prove that the polytope
$D_\sigma$ is always contained in another polytope $P_\sigma$ which is 
the fundamental domain of a certain  affine Weyl subgroup of $\Wa$ associated to  $\g_0$. Hence we can compute 
the volume of
$P_\sigma$ as the index of this subgroup in $\Wa$.\par Finally, it turns out that $P_\sigma\setminus D_\sigma$ is either
empty or it consists of exactly one alcove, so the  volume of $D_\sigma$ can be
computed from that of $P_\sigma$, possibly correcting by $-1$. We finish  the
proof of our enumerative formulas by giving a simple criterion to decide
whenever the correction occurs.
We have 
$$ | \Cal W_{ab}^{\sigma} | = a_0(\chi_\ell(\a_p)+1)k^{n-L}
\frac{|W_f|}{|W_\s|}-\chi_\ell(\a_p),
  $$ where $a_0$ is the label of the vertex $0$ in the Dynkin diagram of type $X_N^{(k)}$,
$\chi_\ell$ is the truth function which is 1 if $\a_p$  is long and 0 otherwise and $L$ is the
number of long simple roots in $\Pi_f$. 
\smallskip
The paper is organized as follows. In section 2 we recall the facts we need about involutions, loop algebras, and Lie
algebra cohomology. In section 3 we introduce the notion of $\s$-minuscule elements and we discuss the relationships of
these elements with
$\frak b_0$-stable abelian subalgebras in $\g_1$. Section 4 is devoted to the study of the polytope $D_\s$. In section 5 and 6 we
prove our main results, in the semisimple and hermitian symmetric case respectively.

\medskip
\heading \S2 Preliminaries
\endheading
\medskip
In this section we summarize the results on affine Kac-Moody  Lie algebras that we shall need in the paper. 
Our main references are \cite{6} and \cite{10}. In particular we use the labeling of Dynkin diagrams as given 
in tables Aff $k,\,k=1,2,3$ in \cite{6, \S4.8}.
\smallskip\noindent{\bf 2.1.}
Given a diagram of type $X_N^{(k)}$, fix a realization $(\ha,\Pia,\Pia^\vee)$  of the
corresponding generalized Cartan matrix. Let $\ga$ denote the corresponding affine Kac-Moody
Lie algebra and $\g$ the finite dimensional Lie algebra of type $X_N$. If $n+2=dim\,\ha$ is
the rank of
$\ga$, then, in
\S~8.6 of
\cite{6}, an automorphism of $\g$ is associated to each $(n+1)$-tuple $\bold
s=(s_0,\dots,s_n)$ of non-negative coprime integers. We denote this
automorphism by $\s_{\bold s,k}$ and call it the automorphism of type $(\bold
s;k)$. The following theorem is Theorem~8.6 of \cite{6}.
\proclaim{Theorem A}
\item a)  $\s_{\bold s;k}$ is of finite order and its order $m$ is given by
$ m=k(\sum_{i=0}^n a_is_i)
$ where $a_i$ are the labels of the diagram $X^{(k)}_N$.

\item b) Up to conjugation by an automorphism of $\g$,  the automorphisms
$\s_{\bold s;k}$ exhaust all $m$-th order automorphisms of
$\g$.

\item c) Two automorphisms $\s_{\bold s;k}$ and $\s_{\bold s';k'}$  are
conjugate by an automorphism of $\g$ if and only if $k=k'$ and the sequence
$\bold s$ can be transformed in the sequence $\bold s'$ by an automorphism of
the diagram $X_N^{(k)}$. 
\endproclaim

Let $\sigma$  be an automorphism of $\g$ of order two and write
$
\g=\g_{0}\oplus\g_{1}
$ for the corresponding gradation.  By the above Theorem, we can assume that
$\s$ is an automorphism of type $(s_0,s_1,\dots,s_n;k)$. Indeed we recover the three
possibilities described in the Introduction. 


Set $\Da$ to denote the set of roots of $\ga$. Let
$\Pia=\{\a_{0},\dots,\a_{n}\}$ be the set of simple roots of
$\ga$ and $\Da^{+}$ the corresponding set of positive roots. If $\a\in\Da$ we
let $\ga_{\a}$ be the corresponding root space.

The numbers $s_{0},\dots,s_{n}$ define a $\ganz$-grading on $\ga$ as follows: if
$\a\in\Da$ write $\a=\sum_{i=0}^{n} m_{i}\a_{i}$ and
$$ ht_{\sigma}(\a)=\sum_{i=0}^{n}s_{i}m_{i}.
$$ Then, if $x\in\ga_{\a}$, we set
$deg(x)=ht_{\sigma}(\a)$. We also set $deg(h)=0$ for $h\in \ha$. We denote by
$
\ga_{i}$ the span of all $x\in\ga$ such that $deg(x)=i$. Set also
$\Da_i=\{\a\in\Da\mid ht_\sigma(\a)=i\}$. Notice that  $\Da_0$
is the root system of $\ga_0$.

\bigskip
\noindent{\bf 2.2.}
Let $L(\g)$ be the loop algebra $\C[t,t^{-1}]\otimes \g$,
$$
\widetilde L(\g)=L(\g)\oplus\C c
$$
its universal (one-dimensional) central extension, and 
$$
\widehat L(\g)=\widetilde L(\g)\oplus \C d
$$
the algebra obtained by extending $\widetilde L(\g)$ by the derivation defined by
$d(p(t)\otimes x)=tp'(t)\otimes x$ and $d(c)=0$.

Let
$\widehat L(\g,\sigma)$ be the subalgebra $\widehat L(\g)$ defined by
$$
\widehat L(\g,\sigma)=\sum_{j\in\ganz}\g_{\overline j}\otimes t^{j}+\C c+\C d
$$ where $\overline j\in\{0,1\}$ is defined by $j\equiv \overline j$ mod $2$.
The following is Theorem~8.5 of \cite{6}.
\proclaim{Theorem B} Let $\h_\s$ be a Cartan subalgebra of $\g_0$. There exists an isomorphism 
$$\Phi:\ga\to\widehat L(\g,\sigma)$$ such that 
\item{i)} $\Phi$ maps $\ga_{i}$ onto $t^i\otimes\g_{\bar
i}$ for $i\ne 0$;
\item{ii)} 
$$
\Phi(\ha)=\h_\s\otimes 1+\C c+\C d;
$$
\item{iii)}
$$
\Phi(\ga_0)=\g_{0}\otimes 1+\C c+\C d.
$$
\endproclaim
\bigskip
\noindent{\bf 2.3.} By means of the grading $ht_\s$, we can  define parabolic subalgebras
$$
\ppa_\s=\ga_0\oplus\ua_\s\text{ and } \ppa_\s^-=\ga_0\oplus\ua_\s^-
$$  of $\ga$, where
$$
\ua_\s=\sum_{ht_\s(\a)>0}\ga_\a\qquad \ua_\s^-=\sum_{ht_\s(\a)<0}\ga_\a.
$$ If we set
$Y=\{i\mid s_i=0\}$, then the parabolic subalgebras $\ppa_\s$ and $\ppa_\s^-$
defined above are the subalgebras
$$
\ppa^{\pm}_Y=\ga_Y\oplus\ua^{\pm}_Y
$$ of $\ga$ as defined in \S~1.2 of \cite{10}.

The grading on $\ga$ defines (by restriction) a  grading on $\ua_\s^-$ and,
henceforth, on
$\L\ua_\s^-$. Let
$(\L^p\ua_\s^-)_q$ denote the subspace of $(\L^p\ua_\s^-)$ of degree $q$. Notice
also that $(\L^p\ua_\s^-)_q=0$ if $q>-p$.

We set $\dd_p:\L^p\ua_\s^-\to\L^{p-1}\ua_\s^-$ to be the standard boundary
operator  affording the Lie algebra homology $H_*(\ua_\s^-)$. It is defined by
setting
$$
\dd_p(x_1\wedge\dots\wedge x_p)=\sum_{i<j}(-1)^{i+j+1}[x_i,x_j]
\wedge\cdots\hat{x}_i\cdots\hat{x}_j\cdots\wedge x_{p}.
$$ if $p>1$ and
$
\dd_1(x)=0
$.
\bigskip We now recall  Garland-Lepowsky generalization \cite{4} of Kostant's theorem. We need
some more notation. Set $\Dap_0=\Dap\cap\widehat \D_0$ and $\Pia_0=\Da_0\cap\Pia$. If $\lambda$ is a dominant weight for
this positive system, denote by 
$V(\lambda)$ be the irreducible $\widehat \g_0$-module of highest weight $\lambda$. Fix $\widehat\rho\in\ha$ satisfying
$\widehat\rho(\a_i)=1$ for all $i,\,0\leq i\leq n$. If $w\in\Wa$ set 
$$N(w)=\{\b\in \Dap\mid w^{-1}(\b)\in -\Dap\}.$$
Recall that  $W_\s$ is the Weyl group of $\widehat \D_0$. Denote by $W_\s'$ the set of elements of minimal
length in the cosets $W_\s w,\,w\in\Wa$. The following is a special case of Theorem 3.2.7 from \cite{10}, which is an extended version of 
Garland-Lepowsky result.
\proclaim{Theorem C}
$$H_p\left(\ua^-_\s\right)=\bigoplus_{{w\in W_\s'}\atop {\ell(w)=p}} V\left(w(\widehat\rho)-\widehat\rho\right)
.$$
Moreover a representative of the highest weight vector of $V\left(w(\widehat\rho)-\widehat\rho\right)$
is given by $e_{-\b_1}\wedge\dots\wedge e_{-\b_p}$ where $N(w)=\{\b_1,\ldots,\b_p\}$  and the $e_{-\beta_i}$ are
root vectors. 
\endproclaim

\bigskip
\noindent{\bf 2.4.} Set 
$$\ha_\R=\sp_\R(\a_1^\vee,\dots,\a_n^\vee)+\R \Phi^{-1}(c)+\R \Phi^{-1}(d).$$
\smallskip\noindent
Define a conjugate linear antihomomorphism  $y\mapsto y^*$
 of $\ga$ setting  $e_i^*=f_i$,
$f_i^*=e_i$ on Chevalley generators  and $h^*=h$ for all $h\in\ha_\R$.  
\smallskip
Let $(\ ,\ )$ be the normalized standard form on $\ga$ \cite{6, \S6.2} and denote by   $\{\ ,\ \}$  
the hermitian form defined by
setting 
$$
\{x,y\}=(x,y^*).
$$  
\smallskip
By \cite{10, Theorem 2.3.13}, the hermitian form  $\{\ ,\ \}$ is positive definite on $\ua^-_\s$. We extend the form 
$\{\ ,\ \}$ to a hermitian form on $\L^p\ua^-_\s$ in the usual way, by determinants.
\bigskip

Let $\dd^*_p:\L^{p-1}\ua_\s^-\to \L^{p}\ua_\s^-$ be the adjoint of $\dd_p$ with respect to the hermitian form $\{\ ,\
\}$  and consider the laplacian $L_p:\L^{p}\ua_\s^-\to \L^{p}\ua_\s^-$
$$ L_p=\dd_{p+1}\dd^*_{p+1}+\dd^*_{p}\dd_{p}.
$$
 Set $\Cal H_p= Ker(L_p)$. Then 
\roster
\item $\Cal H_p\subseteq Ker(\dd_{p})$;
\item the natural map $\Cal H_p\to \left(\Cal H_p\oplus Im(\dd_{p+1})\right)\slash Im(\dd_{p+1})$ induces an isomorphism
$$\Cal H_p\cong H_p(\ua^-_\s).$$
\endroster
\remark{\bf Remark 2.4.1} Notice that $\dd_p(( \L^{p}\ua_\s^-)_q)\subseteq ( \L^{p-1}\ua_\s^-)_q$, and the decomposition
$$\L^{p}\ua_\s^-=\bigoplus_{q\in\ganz} ( \L^{p}\ua_\s^-)_q
$$
is an orthogonal sum. Therefore $\dd_p^*(( \L^{p-1}\ua_\s^-)_q)\subseteq ( \L^{p}\ua_\s^-)_q$. In particular, since
$( \L^{p+1}\ua_\s^-)_{-p}=0$, we have that  $\dd^*_{p+1}=0$ on $( \L^{p}\ua_\s^-)_{-p}$, hence
$${L_p}_{|( \L^{p}\ua_\s^-)_{-p}}={\dd^*_p\dd_p}_{|( \L^{p}\ua_\s^-)_{-p}}.$$
\endremark
\smallskip
\ni{\bf 2.5.} We need some remarks on affine roots.
Recall that  
a root is called real if it is $\Wa$-conjugate to a root in $\widehat\Pi$, imaginary otherwise. 
Imaginary roots are isotropic with respect to  $(\ ,\ )$ whereas real roots are not isotropic. Moreover there are only two 
possible roots 
lengths for real roots, except for the case $\Da\cong A_{2n}^{(2)}$, in which three lengths
occur.  We call long a real root of maximal length and short a real root of minimal length. If only one length occurs, we shall
conventionally say that all roots are long. This convention will be relevant to the
formulation of our results.\par  Set $\d=\sum_{i=0}^na_i\a_i$. The following statements hold.
\roster
\item
 The imaginary roots in $\Da$ are $\pm\Bbb N^+\d$.
\item Suppose that $\ga$ is of type $X^{(k)}_N$. If $\a\in\Da$ then $k\d+\a\in\Dap\cup\{0\}$.
\item If $\a$ is not a long root then $\d+\a\in\Dap\cup\{0\}$.
\endroster
All these properties follow from \cite{6, Theorem 5.6, Proposition 6.3}, where  the
relationships between   the root system $\Da$ and $\D_f$ (the root system generated by
$\Pi_f$) are described in detail. The explicit relation between $\Wa$ and $W_f$ is given in
\cite{6, Proposition 6.5}.
 Sometimes in the following we shall implicitly refer to these descriptions.

\medskip
\heading\S 3 The set $\Cal W^\sigma_{ab}$.\endheading
\medskip

In the following  we identify $\ga$ and $\widehat L(\g,\sigma)$. We observe that,
if $\a\in\Da_0$, then $\a(d)=\a(c)=0$. This  allows us to identify
the set  of roots of
$\g_{0}$ with respect to $\h_{\s}$ with the set $\Da_0$. Recall that 
$\Dap_{0}=\Dap\cap\Da_{0}$ and let
$\bb_{0}$ denote the corresponding Borel subalgebra of $\g_{0}$.
Notice also that
$$
\g_{1}\simeq\ga_{-1}=(\ua_{\s}^-)_{-1}.
$$

\proclaim{Definition} We  say that an element $w\in
\Wa$ is {\sl $\sigma$-minuscule} if
$$ N(w)\subset\{\a\in\Da\mid ht_{\sigma}(\a)=1\}.
$$
\endproclaim
\noindent Denote by $\Cal W^{\sigma}_{ab}$ the set of
$\sigma$-minuscule elements of $\Wa$. We can now state

\proclaim {Theorem 3.1} There is a bijection between  $\Cal W^{\sigma}_{ab}$  and the set $I^\sigma_{ab}$
 of abelian
$\bb_{0}$-stable subalgebras of $\g_{1}$.
\endproclaim
\demo{Proof} Let $\i\subset\g_1$ be a $\bb_0$-stable abelian subalgebra and fix a
basis
$\{x_1,\dots,x_p\}$ of $\i$. Set
$$ v_\i=t^{-1}\otimes x_1\w\dots\w t^{-1}\otimes x_p\in(\L^p\ua_\s^-)_{-p}.
$$

\ni By Remark 2.4.1 we have that
${L_p}_{|(\L^p\hat\u_\s^-)_{-p}}=\dd^*_p\dd_p$. Clearly, since $\i$ is abelian,
$$
\dd^*_p\dd_p(v_\i)=0.
$$ It follows that $v_\i$ is a cycle in $\L^p\ua_\s^-$ and, since $\i$ is
$\bb_0$-stable and $v_\i$ is
$\ha$-stable, its homology class is an highest vector for an irreducible
component $V_\i$  of $H_p(\ua_\s^-)$. By Theorem C there exists an element  $w\in \Wa$ such that $\ell(w)=p$ and
$V_\i=V(w(\hat\rho)-\hat\rho)$. We now check
that $w$ is $\sigma$-minuscule. Suppose that $N(w)=\{\b_1,\dots,\b_p\}$. Then there
is a nonzero $c\in\C$ such that, fixing root vectors $e_{-\b_i}$,
$$ e_{-\b_1}\w\dots\w e_{-\b_p}=c\,\cdot t^{-1}\otimes x_1\w\dots\w t^{-1}\otimes x_p.
$$ Hence $e_{-\b_i}$ lies in the span of the vectors $t^{-1}\otimes x_1,\dots,
t^{-1}\otimes x_p$. This implies that $ht_\sigma(\b_i)=1$.

Thus we have established a map $F:I^\sigma_{ab}\to \Cal W^\sigma_{ab}$. Suppose now
conversely that $w\in \Cal W^\sigma_{ab}$ and set
$$ N(w)=\{\b_1,\dots,\b_p\}.
$$
Since $ht_\sigma(\b_i)=1$ we have that $e_{-\b_i}\in(\ua_\s^-)_{-1}$  hence we can
write
$e_{-\b_i}=t^{-1}\otimes x_i$ with $x_i\in\g_1$. It is well-known that 
$W'_\s=\{w\in\Wa\mid N(w)\cap \Dap_0=\emptyset\}$. In particular, if $w$ is $\s$-minuscule, then 
 $w\in W_\s'$.
Again by Theorem C, the element
$ v=e_{-\b_1}\w\dots\w e_{-\b_p}
$ represents a highest weight vector for $V(w(\hat\rho)-\hat\rho)$ in
$H_p(\ua_\s^-)$. By 2.4 (2) and the subsequent Remark, it follows that
$$ L_p(v)=\dd^*_p\dd_p(v)=0.
$$ It is a standard fact that $\dd^*_p\dd_p(v)=0$ implies $\dd_p(v)=0$. It 
easily follows that the space $\i$ spanned by $\{x_1,\dots, x_p\}$ is  abelian. Since $v$
is
$\bb_0$-stable, then $\i$ is also $\bb_0$-stable.
\endemo

\medskip
\heading \S4 The polytope $D_\sigma$.\endheading
\medskip At this point we wish to count the elements in $I^\sigma_{ab}$ by
counting the elements of $\Cal W^\sigma_{ab}$. This can be done by computing the
volume of certain polytopes.
\smallskip
In the following we identify $\ha_\real$ with $\hard$ via the standard invariant
bilinear form, thus, for all real roots $\a$, $\a^\vee={2\a\over
(\a,\a)}$. Take $\omega_0$ in $\hard$, such that $(\omega_0,\a_0^\vee)=1$, 
$(\omega_0,\a_i^\vee)=0$ for $i\in\{1,\dots,n\}$ and $(\omega_0,\omega_0)=0$.
\bigskip

Set $$\hauno=\{x\in \hard\mid (x,\d)=1\},\qquad \hazero=\{x\in
\hard\mid (x,\d)=0\}.$$ Let $\pi$ be the canonical projection
$\mod \d$ and set $$\huno=\pi \hauno, \qquad\hzero=\pi \hazero.$$
\par
For $x\in \hard$ and $S\subseteq \hard$, we set $\ov x=\pi(x),\,\ov S=\pi(S)$.
We define  a
  $\Wa$-invariant nondegenerate pairing between $\ha_0^*$ and
$\ha^*_\R/\R\d$ by setting
$$ (\a,\l+\R\d)=(\a,\l).
$$
\smallskip
  For
$\a\in\Dap$ set
$$ H_{\a}=\{x\in \h^*_1\mid (\a,x)=0\}
$$ and
$H_\a^+=\{x\in\h_1^*\mid (\a,x)\ge0\}$.
  Set also
$$ C_1=\{x\in \h_1^*\mid (\a,x)\ge0\,\forall\ \a\in\Pia\},
$$ the  fundamental alcove  of $\Wa$.
It is well-known that there is a faithful action of $\Wa$ on $\h_1^*$.
Set
$$ D_\sigma=\bigcup_{w\in\Cal W^\sigma_{ab}} wC_1.
$$ Clearly the number of elements of $\Cal W^\sigma_{ab}$ is equal to
$\frac{\Vol(D_\sigma)}{\Vol(C_1)}$.

Given $w\in \Wa$, a root $\b\in\Dap$ belongs to $N(w)$ if and only if $H_{\b}$ separates $wC_1$ and $C_1$. It follows that
$$ D_\sigma=\bigcap_{{\a\in \hat\Dp,}\atop{ht_\sigma(\a)\ne 1}}H^+_{\a}.
$$
If $\a\in\Da$ we set $\Ha_\a=\{x\in\widehat\h^*_\R/\R\d\mid (\a,x)=0\}$  and
$\Ha^+_\a=\{x\in\widehat\h_\R^*/\R\d\mid (\a,x)\ge0\}$. Set also
$$ C_\sigma=\bigcap_{{\a\in \hat\Dp,}\atop{ht_\sigma(\a)\ne 1}}\Ha_\a^+.
$$ Obviously
$$ D_\sigma=C_\sigma\cap\h_1^*.
$$

For $i\in\{0,\dots,n\}$ we define a number $\epsilon_i$ as follows:
$$\epsilon_i=\cases1\quad&\text{ if $k=1$ or $\a_i$ is not a long root,}\\
2&\text{ if $k=2$ and $\a_i$ is a long root.}\endcases
$$

\ni Denote by
$\Da_0^{max}$  the set of maximal roots in $\Dap_0$. Set
$$\Phi_\sigma=\{\a_i+\epsilon_i s_i\d\mid
i=0,\dots,n\}\cup\{k\d-\gamma\mid\gamma\in\Da_0^{max}\}.
$$
\proclaim{Proposition 4.1} We have  $\Phi_\sigma\subset\Dap$ and
$$ C_\sigma=\bigcap_{\a\in\Phi_\sigma}\Ha^+_\a.
$$
\endproclaim
\demo{Proof} Set
$$ P=\bigcap_{\a\in\Phi_\sigma}\Ha^+_\a.
$$ First of all we prove that
$\Phi_\sigma\subseteq\{\a\in\Dap\mid ht_\sigma(\a)\ne1\}$. This will imply obviously
that $C_\sigma\subseteq P$.

It is clear that, if $s_i=0$, then $\a_i=\a_i+\epsilon_is_i\d\in\Dap$ and
$ht_\sigma(\a_i)=0\ne1$.  By 2.5 we have that
$\Phi_\sigma\subseteq\Da\cup\{0\}$. We observe that
$$ht_\sigma(\d)=ht_\sigma(\sum_{i=0}^n a_i\a_i)=\sum_{i=0}^n a_is_i=\frac{2}{k},$$ hence, if
$s_i\ne0$,
$ht_\sigma(\a_i+\epsilon_is_i\d)=s_i(\frac{2}{k}\epsilon_i+1)>1$. It follows
that $\a_i+\epsilon_is_i\d\in\Dap$ and $ht_\sigma(\a_i+\epsilon_is_i\d)\ne1$. In
the same manner, if $\gamma\in\Da_0^{max}$, then
$ht_\sigma(k\delta-\gamma)=2>1$, hence $k\delta-\gamma\in\Dap$ and
$ht_\sigma(k\d-\gamma)\ne1$.

We now check that $P\subseteq C_\sigma$.  If $x\in P$ and $\a\in\Dap$, then,
writing $\a=\sum_{i=0}^n n_i\a_i$, we obtain 
$$ 0\le\sum_{i=0}^n n_i(\a_i+\epsilon_is_i\d,x)=(\a,x)+(\sum_{i=0}^n n_i\epsilon_is_i)(\d,x).
$$ If $\a=\d$, then $(1+\sum_{i=0}^n a_i\epsilon_is_i)(\d,x)\ge0$ hence
$(\d,x)\ge0$. If $k=1$  then we can rewrite the above formula as
$$ (\a,x)\ge-(\sum_{i=0}^n n_is_i)(\d,x)=-ht_\sigma(\a)(\d,x).
$$ If $k=2$ then there is a unique index $p$ such that $s_p\ne0$, so we can write
$$ (\a,x)\ge-\epsilon_p(\sum_{i=0}^n n_is_i)(\d,x)\ge-2ht_\sigma(\a)(\d,x).
$$ In any case we have
$$ (\a,x)\ge-k\,ht_\sigma(\a)(\d,x).
$$ In particular, if
$\a\in\Dap_0$, then
$(\a,x)\ge0$.

If $\a\in\Da$ is a real  root then there is $m\in\ganz$ such that
$\a=km\d+\b$, with $\b\in\Da_0\cup\Da_1$. Indeed, if $ht_\sigma(\a)=r$, we can
write $\a=k[r/2]\d+\a-k[r/2]\d$, so we can choose $m=[r/2]$ and
$\b=\a-k[r/2]\d$. Notice also that, if $\a\in\Dap$, then $m\ge0$. Therefore,  if $\a\in\Dap$ is such that $ht_\sigma(\a)\ne1$,
there are the following four possibilities.

\item{A)}$\a=m\d$ with $m>0$. If $x\in P$  then we have seen above that
$(\d,x)\ge0$.

\item{B)}$\a=km\d+\b$ with $\b\in\Dap_0$ and $m\ge0$. If $x\in P$ then
$(\a,x)=km(\d,x)+(\b,x)\ge0$.

\item{C)}$\a=km\d-\b$ with $\b\in\Dap_0$ and $m\ge1$. There is a root
$\gamma\in\Da_0^{max}$ such that $\gamma-\b$ is a sum of roots in $\Dap_0$. This
implies that, if $x\in P$, then
$(\a,x)=km(\d,x)-(\b,x)=k(m-1)(\d,x)+(\gamma,x)-(\b,x)+k(\d,x)-(\gamma,x)\ge0$.

\item{D)}$\a=km\d+\b$ with $\b\in\Da_1$ and $m\ge1$. If $x\in P$ then
$km(\d,x)+(\b,x)\ge(km-k)(\d,x)\ge0$.
\endemo
\bigskip

Let $\G_\s$ be the Dynkin graph of $\Pia_0$. $\G_\s$ is, in
general, disconnected. Each of its connected components is of
finite type. We write $\Si|\G_\si$ if $\Si$ is a connected
component of $\G_\si$. Assume that $\Si|\G_\si$, and denote by
$\Pi_\Si$ the simple roots in $\Si$, by $W_\Si$ the relative Weyl
group, $W_\Si=\langle s_\a\mid \a\in\Pi_\Si\rangle$, and  by
$\D_\Si$ the relative root system, $\D_\Si=W_\Si\Pi_\Si$.
Moreover, let $\th_\Si$ be the highest root of $\D_\Si$, 
$\a_\Si=k\d-\th_\Si$, $\Pia_\Si=\Pi_\Si\cup\{\a_\Si\}$, and
$\What_\Si=\langle s_\a\mid \a\in\Pia_\Si\rangle$. If $X_{|\Si|}$ is the type of the (finite) system $\D_\Si$, then the root system
generated by $\Pia_\Si$ is clearly of type $X_{|\Si|}^{(1)}$.
We set $\Pia_\si=\bigcup\limits_{\Si|\G_\si}\Pia_\Si$. 
Since $\Da^{\text{\it max}}_0=\{\th_\Si\mid \Si|\G_\s\}$, we have
$$\Phi_\s=\Pia_\s\cup\{\a_i+\epsilon_i s_i\d\mid \a_i\notin\Pia_0\}.$$
We also set 
$$W_\si=\langle s_\a\mid \a\in\Pia_0\rangle, \qquad
\What_\si=\langle s_\a\mid \a\in\Pia_\si\rangle.$$ $W_\si$
and $\What_\si$  are the direct products
$\prod\limits_{\Si|\G_\si}W_\Si$ and
$\prod\limits_{\Si|\G_\si}\What_\Si$.
\smallskip
Observe  that $rk\,\g_0=rk\,\ga-2$ while the rank of
$[\g_0,\g_0]$ is given by the rank of the subsystem $\Da_0$. It
follows that $\g_0$ is semisimple if and only if there is a
simple root $\a_p\in\Pia$ such that
$\Pia_0=\Pia\backslash\{\a_p\}$. The calculation of the order of
$\Cal W_{ab}^\sigma$ is better performed by separating the case
when $\g_0$ is semisimple from the case when $\g_0$ has  a
nontrivial center.
\par 
To simplify notation, henceforward we identify $\hzero$ with
$\sum\limits_{j=1}^n\real \a_j$ and $\pi(\omega_0)$ with
$\omega_0$. Hence $\huno=\omega_0+\hzero$. As usual, we set
$\th=\d-a_0\a_0$.

\medskip
\heading \S5 The semisimple case\endheading
\medskip
We assume that $\Pia_0=\Pia\backslash\{\a_p\}$, hence $\Phi_\s
=\Pia_\s\cup\{\a_p+s_p\epsilon_p\d\}$. We define
$$P_\s=\bigcap\limits_{\a\in\Pia_\s} H_\a^+.$$ 
Then $D_\s\subseteq P_\s$. We shall first compute
$\Vol(P_\s)/\Vol(C_1)$; then we shall see that the set
difference $P_\s\setminus D_\s$ is either  empty or 
exactly one alcove of $\Wa$. It is easily seen that $P_\s$ is a
fundamental domain for $\Wa_\s$, hence by \cite{1, VI.4, Lemma 1}
we have that $\Vol(P_\s)/\Vol(C_1)= [\Wa:\Wa_\s]$.  This
description is not directly useful to obtain an explicit result.
Indeed, we shall compute the volume of a certain {\it translate}
of $P_\s$. We need some preliminaries about translations and 
the structure theory of $\Wa$.
\par
For $\a\in\hzero$, we define the {\it translation} by $\a$,
$t_\a:\hard\to \hard$, setting
$$
t_\a(x)=x+(x,\d)\a-((x,\a)+\frac{1}{2}|\a|^2(x,\d))\d.
$$ 
This agrees with Kac's definition \cite{6, 6.5.2}, according to
our conventions. The translation $t_\a$ preserves the standard
invariant form for all $\a\in \hzero$. Since $t_\a(\d)=\d$ for
all $\a\in\hzero$, $t_\a$ induces a map $\pi(\hard)\to
\pi(\hard)$ which we still denote by $t_\a$. Notice that for
$x\in\h^*_1$ and $y\in \hazero$ we have
$$
t_\a(x)=x+\a,\qquad t_\a(y)=y-(y,\a)\d.
$$
For any $S\subseteq \hzero$, we set $T(S)=\{t_\a\mid\a\in S\}$. 
Then $\What=W_f\ltimes T(M)$, with $M$ the lattice generated
by $W_f(\th^\vee)$ \cite{6, Proposition 6.5}. We explicitly 
describe the lattice $M$. The first part of the next proposition
is well known.

\proclaim{Proposition 5.1} {\rm (1)}\ If $\Da$ is untwisted, or
$a_0=2$, then $M=Q^\vee$. 
\par\ni
{\rm (2)}\ If $k>1$ and $a_0=1$, then
$$
M=\sum\limits_{i=1}^n r_i\ganz\a_i^\vee,\eqno 
$$
with $r_i=1$ if $\a_i$ is short, and $r_i=k$ if $\a_i$ is long.
\endproclaim

\demo{Proof}
(1) If $\Da$ is untwisted, or $a_0=2$, then $\th$ is long, so that
$W_f(\th^\vee)$ is the set of short coroots. Since $Q^\vee$ is
generated by the short coroots, we obtain  $M=Q^\vee$.  
\par
\ni
(2) If $k>1$ and $a_0=1$, then $\th$ is a short root and since
$w(\th^\vee)=w(\th)^\vee$ for all $w\in W_f$,  $M$ includes
$\a^\vee$ for all short $\a\in\D_f$. Any long $\b\in \D_f$ is an
integral linear combination of short roots, say
$\b=c_1\b_1+\cdots+c_s\b_s$, and $\b^\vee={2\over(\b,\b)}
(c_1\b_1+\cdots+c_s\b_s)={(\th,\th)\over(\b,\b)}
(c_1\b_1^\vee+\cdots+c_s\b_s^\vee)$. Since for any long $\b$ we have
that ${(\b, \b)\over(\th, \th)}=k$, we obtain  $k\b^\vee\in M$ 
for all long $\b$.
Set $M'=\sum\limits_{i=1}^n r_i\ganz\a_i^\vee,$ with $r_i=1$ if
$\a_i$ is short, and $r_i=k$ if $\a_i$ is long. Then we have that
$M'\subseteq M$. We shall prove that in fact $M=M'$. It suffices 
to prove that for any short $\a_i$ and $w\in W_f$ we have 
$w(\a_i^\vee)\in M'$, since $W_f\a_i^\vee$ generates $M$. 
It is 
clear that if $\a_j$ and $\a_i$ are short roots, then $s_{\a_j}
(\a_i^\vee)\in M'$. If  $\a_j$ is long and $\a_i$ is short, 
then $s_{\a_j}(\a_i^\vee)=\a_i^\vee-(\a_i^\vee, \a_j)\a_j^\vee$ 
and since $k| (\a_i^\vee, \a_j)$, $s_{\a_j}(\a_i^\vee)\in M'$. 
Since moreover $W_f(k \a^\vee)\subseteq M'$ for all $\a\in\D_f$,
we inductively obtain  $W_f\a_i^\vee\subseteq M'$.
\endemo
\bigskip

The group $\Wa_\s$ is not in general a product of a subgroup of $W_f$ and 
a subgroup of $T(M)$. This motivates the following 
construction. 
\par

Assume that $\Si|\G_\si$. If $\a_0\in\Pi_\Si$, we set
$\Pi_{\Si,f}=\Pi_\Si\setminus\{\a_0\}\cup \{-\th\}$, and if
$\a_0\notin\Pi_\Si$, we set $\Pi_{\Si, f}=\Pi_\Si$. Moreover, we
set $\Pia_{\Si,f}=\Pi_{\Si,f}\cup \{k\d-\ov\th_\Si\}$ and denote
by $\D_{\Si,f}$, $\Da_{\Si,f}$ the root systems generated by
$\Pi_{\Si,f}$, $\Pia_{\Si,f}$, respectively.
\par

If $\Da\not\cong A^{(2)}_{2n}$, it is clear that $\ov \th_\Si\in
\D_f$ and $k\d-\ov\th_\Si\in \Da$ (see \cite{6}, warning after
6.3.8). Moreover, $\ov\th_\Si$ is the highest root of
$\D_{\Si,f}$. Therefore $\D_{\Si,f}$ and $\Da_{\Si,f}$ are
isomorphic to $\D_\Si$ and $\Da_\Si$, respectively. Moreover
$\Da_{\Si,f}$ is isomorphic to the untwisted affine system
associated to $\D_{\Si,f}$.
\par

Assume that $\Da\cong A^{(2)}_{2n}$. Then $p=n$, and there is
just one connected component $\Si=\G_\si$. $\D_{\Si}$ is of type
$B_n$, but $\D_{\Si,f}$ is of type $C_n$. In fact, $\Pi_{\Si,f}=
\{-\th, \a_1, \dots, \a_{n-1}\}$, $\a_1, \dots,\a_{n-1}$ are 
short roots of $\D_{\Si,f}$, and $-\th$ is long.  Moreover,
$\th_\Si=2\a_0+\cdots+2\a_{n-2}+\a_{n-1}$, and hence $\ov
\th_\Si= -\th+2\a_1+\cdots+2\a_{n-2}+\a_{n-1}$. Thus $\ov\th_\Si$
is the highest short root of $\D_{\Si,f}$ and $\Da_{\Si,f}$ is
a twisted affine root system of type $A^{(2)}_{2n-1}$.
\par

We define $W_{\Si,f}=\langle s_\a\mid \a\in \Pi_{\Si,f}\rangle$
and $\What_{\Si, f}=\langle s_\a\mid\a\in \Pia_{\Si, f} \rangle$.
Since $s_{-\th}=s_{\ov\a_0}$, $W_{\Si,f}=\langle s_\a\mid \a\in 
\ov \Pi_{\Si}\rangle$. Moreover, since for any non-isotropic
$\a, \b\in \hazero$ $s_\a s_\b$ and $s_{\ov \a} s_{\ov\b}$, have
the same period, we obtain that $\What_\Si$ and $\What_{\Si,f}$,
with the given sets of generators, are naturally isomorphic as
Coxeter systems, though the respective root systems are not
necessarily isomorphic.
\par

Finally, we set
$\Pi_{\si,f}=\bigcup\limits_{\Si|\G_\si}\Pi_{\Si,f}$,
$\Pia_{\si,f}=\bigcup\limits_{\Si|\G_\si}\Pia_{\Si,f}$, and
$$W_{\si,f}=\langle s_\a\mid \a\in\Pi_{\si,f}\rangle, \qquad
\What_{\si,f}=\langle s_\a\mid \a\in\Pia_{\si,f}\rangle.$$
$W_{\si,f}$ is the direct product of the (finite) Weyl groups
$W_{\Si,f}$, and $\What_{\si,f}$ is  the direct product of the
$\What_{\Si,f}$, for all $\Si|\G_\s$.
\bigskip

\proclaim{Proposition 5.2} We have
$\What_{\si,f}=W_{\si,f}\ltimes T(M_\si),$
where
$$M_\si=\cases \sum\limits_{\a\in\Pi_{\si,f}} k\ganz\a^\vee 
\quad &\text{ if } \quad \Da\not\cong A^{(2)}_{2n}\\
4\ganz\th^\vee\oplus 2\ganz\a_1^\vee\oplus\cdots\oplus
2\ganz\a_{n-1}^\vee \quad &\text{ if } \quad \Da\cong
A^{(2)}_{2n}\\ \endcases$$
\endproclaim

\demo{Proof}
By \cite{6, Proposition 6.5} $\Wa_{\Si,f}=W_{\Si,f}\ltimes
T(M_\Si)$, where $M_\Si$ is the lattice generated by
$W_{\Si,f}(k\ov \th^\vee_\Si)$. Hence the claim follows by the above
discussion and Proposition 5.1.
\endemo
\bigskip

The above structure results allows us to compute explicitly  the
index of $\Wa_{\s, f}$ in $\Wa$.

\proclaim{Proposition 5.3} We have
$$[\What:\What_{\si,f}]= r_p k^{n-L} \ [W_f:W_{\si,f}],$$ 
where $L$ is the number of long roots in $\Pi_f$ and $r_p$ is the
ratio between the squared length of $\a_p$ and that of any short
root.
\endproclaim

\demo{Proof} 
By standard group theory we obtain  $[\What:\What_{\si,f}]=
[W_f:W_{\si,f}][M:M_\si]$. We have only to prove that, in all 
cases, $[M:M_\si]=r_p k^{n-L}$.
\par

First assume that  $p=0$. This implies that  $k=2,\,a_0=1$. Then
$W_f=W_{\s,f}$ and $M_\sigma=\sum\limits_{\a\in\Pi_f}k\ganz
\a^\vee$. Moreover, $r_p=1$. Hence the claim follows by 
Proposition 5.1~(2).
\par 

Next assume that $p\geq 1$. Set $\th^\vee=b_1\a_1^\vee+\dots
+b_n\a_n^\vee$. Suppose that $k=1$, so that $a_p=2$. If $\a_p$
is short, then $b_p=1$, hence $\a_p^\vee\in
\sum\limits_{\a\in\Pi_{\si,f}} \ganz\a^\vee=M_\si$. It follows
that $M_\si=Q^\vee=M$, hence the claimed equality holds. If
$\a_p$ is long, then $b_p=2$ and $2\a_p^\vee \in M_\si$, but
$\a_p^\vee\notin M_\si$. It follows that $[Q^\vee:M_\si]=2$,
hence the claim is true in this case, too.
\par

Now we assume  that $k=2$, so that $a_p=1$. If $\a_p$ is short, then
$b_p=1$. Moreover, $\Da\not\cong A^{(2)}_{2n}$. Thus we obtain
 $2\a_p^\vee\in M_\si$. It follows that $M_\si=2 Q^\vee$, hence
$[M:M_\si]=2^R$, where $R$ is the number of short roots in $\Pi_f$,
and the claim holds. Then we assume that $\a_p$ is long.  
If $\Da\not\cong A^{(2)}_{2n}$, then $b_p=2$, hence
$4\a_p^\vee\in M_\si$ and $2\a_p^\vee\notin M_\si$. It follows that
$M_\si=4\ganz\a_p^\vee+\sum\limits_{{\a\in\Pi_f}\atop{\a\ne\a_p}}
2\ganz\a^\vee$, hence $[M:M_\si]=2^{R+1}$, where $R$ is the number of
short roots in $\Pi_f$, which is equivalent to our claim.
Finally, if $\Da\cong A^{(2)}_{2n}$, then $p=n$ and $\a_n$ is
long. We have $b_n=1$, hence $4\a_n^\vee\in M_\si$ and $2
\a_n^\vee\notin M_\si$. It follows that $M_\si=2\ganz\a_1^\vee
\oplus \cdots \oplus 2\ganz\a_{n-1}^\vee \oplus 4\ganz\a_n^\vee$,
hence $[M:M_\si]=2^{n+1}$, which is the claim in this case.
\endemo
\bigskip

We then prove that $\Vol(P_\s)$ is equal to $[\Wa:\Wa_{\s,f}]$.
Set $$A_{\s,f}=\bigcap\limits_{\a\in\Pia_{\s,f}}H_\a^+.$$
It is clear that $H_{-\th}^+=H_{\ov \a_0}^+$, hence
$A_{\s,f}=(\bigcap\limits_{i\ne p}H_{\ov\a_i}^+)\cap 
H^+_{k\d-\ov \th_\Si}$. 

\proclaim{Lemma 5.4}
$A_{\s,f}$ is a fundamental domain for the action of
$\What_{\si,f}$ on $\huno$. Hence
$$
{\Vol\left(A_{\s,f}\right)\over\Vol(C_1)}=[\Wa:\Wa_{\s,f}].
$$
\endproclaim

\demo{Proof}
For any $\Si|\G_\si$, set $\h^*_{\Si,0} =\sum\limits_
{\a\in\Pi_{\Si,f}}\real \a$ and $A_{\Si,f}=\{x\in
\h^*_{\Si,0}\mid (x,\a)\geq 0\ \forall\, \a\in\Pi_{\Si,f}; \
(x,\ov\th_\Si)\leq k\}$. Then we have an orthogonal decomposition
$\hzero=\sum\limits_{\Si|\G_\si} \h^*_{\Si,0}$, and since
$\huno=\omega_0+\hzero$ we obtain $A_{\s,f}=\omega_0+
\sum\limits_{\Si|\G_\si}A_{\Si,f}$. Now $\What_{\Si,f}$ acts
faithfully onto $\omega_0+ \h^*_{\Si,0}$, and $\omega_0+
A_{\Si,f}$ is its fundamental alcove, hence a fundamental domain
for this action. Moreover $\What_{\Si, f}$ fixes pointwise
$\h^*_{\Si',0}$, for all other $\Si'|\G_\s$. Since $\What_{\si,f}$
is the direct product $\prod\limits_ {\Si|\G_\si} \What_{\Si,
f}$, we obtain that $\omega_0+\sum \limits_{\Si|\G_\si}
A_{\Si,f}$ is a fundamental domain for the action of $\What_{\si,f}$
on $\huno$.
\endemo
\medskip

Let $\{\omega_j^\vee\mid 1\leq j\leq n\}$ be the dual basis of
$\Pi_f$ in $\hzero$ and $o_j=\displaystyle{\omega_j^\vee\over
a_j}$ for $1\leq j\leq n$, where the $a_i$ are the labels of
the Dynkin diagram. Thus we have $(\th, o_j)=1$.
Set moreover $o_0=0$ and for $j\in \{0,1,\dots,n\}\setminus
\{p\}$ define $\tilde\omega_j^\vee= a_j(o_j-o_p)$. Then $\{ \tilde
\omega^ \vee_j \mid j\in \{0,1,\dots,n\}\setminus \{p\}\}$ is the
dual basis of $\{\ov\a\mid\a\in\Pia_0\}$ in $\hzero$. In particular
it is a basis of $\hzero$. Indeed, $\{ \tilde \omega^ \vee_j \mid
\a_j\in\Pi_\Si\}$ is the dual basis of $\ov \Pi_{\Si}$ in 
$\h^*_{\Si,0}$, for all $\Si|\G_\si$. 
\par
For all $\Si|\G_\si$, let $\th_\Si=\sum\limits_{\a_j\in
\Pi_\Si} a_j'\a_j$: this defines integers $a_j'$ for all $j\in
\{0,1,\dots,n\}\setminus\{p\}$. We set $\tilde o_j=
k\displaystyle {\tilde\omega^\vee_j\over a'_j}$, for all $j\in
\{0,1,\dots,n\}\setminus \{p\}$, so that $(\tilde o_j, \th_\Si)=
(\tilde o_j, \ov\th_\Si)=1$, for all $j$ such that $\a_j\in \Pi_\Si$.
\bigskip

\proclaim{Proposition 5.5} We have
$P_\s=o_p+ A_{\s,f}$, hence
$$\Vol(P_\s)=[\Wa:\Wa_{\s,f}].$$ 
Moreover, $t_{-o_p}\What_\si t_{o_p}= \What_{\si,f}.$
\endproclaim

\demo{Proof} We have 
$$ t_{-o_p}(\a_j)=\a_j+(o_p, \a_j)\d=
\cases \a_j\qquad &\text{if $ j\notin\{0,p\}$,}\\
       \a_0-{\d\over a_0}=-{\th\over a_0} \qquad &\text{if $j=0$.}
\endcases$$
Thus, in any case, $t_{-o_p}(\a_j)=\ov \a_j$, and hence
$t_{-o_p}(k\d-\th_\Si)= k\d-\ov\th_\Si$. Since $t_{-o_p}$
preserves the standard invariant form, this implies 
that $t_{-o_p}(H^+_{\a_j})=H^+_{\ov \a_j}$ for all $j\ne p$, and 
$t_{-o_p}(H^+_{k\d-\th_\Si})=H^+_{k\d-\ov \th_\Si}$, for all 
$\Si|\G_\s$. By the definitions of $P_\s$ and $A_{\s, f}$, it
follows that $t_{-o_p}(P_\s)=A_{\s,f}$, hence the claim.
\par

Since $t_{-o_p}$ preserves the standard invariant form, we also
have that $t_{-o_p} s_\a t_{o_p}$ $= s_{t_{-o_p}(\a)}$. Since
$s_\th=s_{\th\over a_0}$, we obtain that conjugation by $t_{o_p}$
maps the generators of $\What_\si$ onto those of $\What_{\si,f}$, 
hence it maps $\What_\si$ onto $\What_{\si,f}$.
\endemo
\bigskip

We are now going to study the difference set $D_\s\setminus
P_\s$.

\remark{\bf  Remark } 
(1). Notice that $k=\frac{2}{a_p}$.
\smallskip\ni
(2). If $\a_p$ is long, then $\epsilon_p s_p=k$, hence
$D_\s=P_\s\setminus H_{\a_p+k\d}^-.$
\bigskip
\endremark 

\proclaim{Lemma 5.6}
Assume that $\a_p$ is long. Then, for each $\Si|\G_\si$, $\a_p$ is
connected to exactly one root $\a_{j(\Si)}$ in $\Si$. Moreover,
$a'_{j(\Si)}=1$, and $-k\a_p^\vee=\sum\limits_{\Si|\G_\si} \tilde
o_{j(\Si)}$.
\endproclaim

\demo{Proof}
Let $\Si|\G_\si$. Our assumptions imply that $\Da\not\cong
A^{(1)}_n$, $n\geq 2$, hence $\a_p$ is connected to exactly one
root in $\Si$, say $\a_{j(\Si)}$. Since $\a_{j(\Si)}$ supports
$\th_\Si$ and any other root in $\supp\th_\Si$ is orthogonal to
$\a_p^\vee$, we also obtain $(\a_p^\vee,\th_\Si)\ne 0$.
Since $\a_p$ is a long root, we have that
$(\a_p^\vee,\a_{j(\Si)})=-1$ and $|(\a_p^\vee,\th_\Si)|\leq 1$;
hence $(\a_p^\vee,\th_\Si)=-1$ and therefore $a'_{j(\Si)}=1$. It
follows that $(\tilde o_{j(\Si)}, \a_{j(\Si)})=k$ and that
$-k\a_p^\vee=\sum\limits_{\Si|\G_\si} \tilde o_{j(\Si)}$.
\endemo
\bigskip

The above results imply in particular that if $\a_p$ is long,
then $D_\s$ is properly included in $P_\s$, since $-k\a_p^\vee\in
P_\s$ and $(-k\a_p^\vee, \a_p)=-2k<-k$.
\smallskip

\proclaim{Proposition 5.7} We have
$D_\s=P_\s$ if and only if $\a_p$ is short. 
\endproclaim

\demo{Proof} 
It suffices to  prove that if $\a_p$ is short then $D_\s=P_\s$.
Since $\a_p$ is short, for at least one $\Si|\G_\si$ we have that
$\th_\Si$ is a long root, and therefore $(\th_\Si,
\a_p^\vee)<-1$. It follows that $\th_\Si+2\a_p$ is a root. Hence
$k\d-\th_\Si-2\a_p\in \Da_0$. Now we have $2\a_p-k\d=-(k\d-
\th_\Si-2\a_p)-\th_\Si$ and since $(x, \b)< k$ for all positive
roots $\b\in \Da_0$ and $x\in P_\s$, we obtain
$(x,2\a_p-k\d)>-2k$. It follows that  $(x, \a_p)\geq-{1\over 2}k$
for all $x\in P_\s$, hence $D_\s=P_\s$.
\endemo
\bigskip

Finally we show that when $D_\s\ne P_\s$ then the difference set
is an alcove of $\Wa$. In fact, we shall explicitly provide an
element $w_\s\in\Wa$ such that $w_\s(C_1)=P_\s\setminus D_\s$.
\par

For any $\Si| \Gamma_\Si$ set
$$t_\Si=t_{\tilde o_{j(\Si)}},\qquad w_\Si=w_\Si^*w_\Si^0$$ 
where  $w_\Si^0$ is the longest element in $W_\Si$ and $w_\Si^*$
is the longest element in
$W(\Pi_\Si \setminus \{\a_{j(\Si)}\})$. Define then
$$w_\si=\prod_{\Si| \Gamma_\Si}t_\Si w_\Si.$$

Next lemma implies in particular that $w_\s\in \Wa$.

\proclaim{Lemma 5.8} We have
$$w_\s(P_\s)=P_\s.$$
Moreover 
$$w_\s=t_{-k\a_p^\vee} w_0^*w_0$$ 
where $w_0$ is the longest element in $W_\s$ and $w_0^*$ is the
longest element of
\break  
$W(\Pia\cap\{\a_p\}^\perp)$.
\endproclaim

\demo{Proof} 
By Lemma 5.6 $a'_{j(\Si)}=1$, hence (see \cite{5, Section 1}) 
$w_\Si(\Pi_\Si\cup\{-\th_\Si\})=\Pi_\Si\cup\{-\th_\Si\}$.
Moreover, $w_\Si(-\th_\Si)=\a_{j(\Si)}$. Since $\tilde
o_{j(\Si)}\perp\{\a\in\Pi_\Si\mid\a\ne\a_{j(\Si)}\}$ and $(\tilde
o_{j(\Si)},\a_{j(\Si)})=(\tilde o_{j(\Si)},\th_\Si)=k$, it
follows that $t_\Si w_\Si(\Pia_\Si)=\Pia_\Si$. Moreover, $t_\Si
w_\Si$ fixes pointwise $\Pia_{\Si'}$ for $\Si'\,| \G_\s,
\,\Si'\ne\Si$, hence $t_\Si w_\Si(\Pia_\s)=\Pia_\s$. It follows
that $w_\si(\Pia_\s)=\Pia_\s$, and hence that $w_\s(P_\s)=P_\s$.
\par
Since $w_\Si t_{\Si'}=t_{\Si'}w_\Si$ for all $\Si, \Si'|\G_\s$ 
with $\Si\ne \Si'$, it is clear that $w_\s=\prod_{\Si| \Gamma_\Si}
t_\Si w_0^* w_0$, hence by Lemma 5.6 $w_\s=t_{-k\a_p^\vee}
w_0^*w_0$.
\endemo
\bigskip

It is well known  that $-w_0$ induces a permutation of $\Pia_0$.
Since $\Pia_0=\Pia\backslash\{\a_p\}$ we can define
$\a_{i'}=-w_0(\a_i)$ for $i\ne p$. For calculating the action of
$w_\s$ on $C_1$ we need the following lemma.

\proclaim{Lemma 5.9}
$$a_i=a_{i'}.$$
\endproclaim
\demo{Proof} 
We first prove that $\g_1$ is irreducible as $\g_0$-module. 
Remark that, as  $\g_0$-modules, $\g_1\cong (\L^1\ua_\s^-)_{-1}$.
Since $L_1=\dd_1^*\dd_1=0$, we see that
$(\L^1\ua_\s^-)_{-1}$ is a nontrivial submodule of $H_1(\ua_\s^-)$, but, 
by Theorem C,
$$
H_1(\ua_\s^-)=\bigoplus_{{w\in
W'_\Si}\atop{\ell(w)=1}} V(w(\widehat\rho)-\widehat\rho) =
V(s_{\a_p}(\widehat\rho)-\widehat\rho))=V(-\a_p).
$$
Moreover, the highest  weight of $\g_1$ as a
$\g_0$-module is $-\a_p$ restricted to $\sp(\a_i^\vee\mid i\ne
p)$, which is equal to $\frac{1}{a_p}\sum_{i\ne p}a_i\a_i$. Since
$\g^*_1\cong \g_1$ as $\g_0$-modules, applying $w_0$ gives the
desired result.
\endemo 
\bigskip

\proclaim{Lemma 5.10} We have
$$w_0t_{k\a^\vee_p}(\a_p)=-k\d-\a_p$$
hence, in particular, $w_\s(C_1)=P_\s\backslash D_\s$ and
$$
\Vol(P_\s\backslash D_\s)=\Vol(C_1).
$$
\endproclaim

\demo{Proof} 
We have $t_{k\a_p^\vee}(\a_p)=\a_p-(k\a_p^\vee,
\a_p)\d=\a_p-2k\d$. To prove the Lemma we have to check that
$w_0(\a_p)=k\d-\a_p$. Let $\{\widehat\omega_0, \dots,
\widehat\omega_n,\d\}$ be the dual basis of $\{\a_0, \dots, \a_n,
\omega_0\}$.  It suffices to check that
$(w_0(\a_p),\widehat\omega_i)=ka_i-\d_{ip}$ ($\d_{ip}$ is a
Kronecker $\d$). For $i=p$ this is obvious. Assume that $i\ne p$. We
have $(w_0(\a_j),\widehat\omega_i)=(\a_j,w_0(\widehat\omega_i))$
for $0\leq j\leq n$. This implies that
$w_0(\widehat \omega_i)=-\widehat\omega_{i'}+ m\widehat\omega_p+
r\delta$ for some $m,r \in \ganz$. By applying $\d$ to both sides
of the previous equation we find that $a_i=-a_{i'}+a_pm$, hence, by Lemma 5.9, that 
$m=\frac{2a_i}{a_p}=k a_i$. But clearly
$(w_0(\a_p),\widehat\omega_i)=m$, hence we get the claim.
\endemo
\bigskip

Putting together  5.3,~5.5, and 5.10  we obtain  the main
result of this section.

\proclaim{Theorem 5.11} 
Assume that $\g_0$ is semisimple. If $\chi_\ell(\a_p)$ is the
truth function which is 1 if $\a_p$ is long and 0 otherwise, then
$$ 
|\Cal  W_{ab}^{\sigma} | = a_0(\chi_\ell(\a_p)+1)k^{n-L}
\frac{|W_f|}{|W_\s|}-\chi_\ell(\a_p)
$$
 where $L$ is the number of long roots in $\Pi_f$.
\endproclaim
\bigskip


The uniform formula established in the previous theorem can be made completely 
explicit in each case. If $k=1$ and $a_p=2$, we have

$$
\alignat5 &\text{type of $\widehat\g$}&&\qquad p &&\text{type of $\D_f$}\ &&\
\text{type of
$\g_0$}&&|\Cal W_{ab}^{\sigma} |\\ &B_n^{(1)}\qquad\quad&& 2\leq p\leq
n-1\quad\qquad&&B_n\quad\qquad&&D_p\times B_{n-p}\quad\qquad&&4\binom{n}{p}-1\\
&\qquad&&p=n\qquad&&B_n\qquad&&D_n\qquad&&2\\ &C_n^{(1)} \qquad&& 1\leq p\leq
n-1\qquad&&C_n\qquad&&C_p\times C_{n-p}\qquad&&\binom{n}{p}\\ &D_n^{(1)}
\qquad&& 2\leq p\leq n-2\qquad&&D_n\qquad&&D_p\times
D_{n-p}\qquad&&4\binom{n}{p}-1\\ &G_2^{(1)}\qquad&&
p=1\qquad&&G_2\qquad&&A_1\times A_1\qquad&&5\\ &F_4^{(1)}\qquad&&
p=1\qquad&&F_4\qquad&&A_1\times C_3\qquad&&23\\ &\qquad&&
p=4\qquad&&F_4\qquad&&B_4\qquad&&3\\ &E_6^{(1)}\qquad&&
p=2,4,6\qquad&&E_6\qquad&&A_1\times A_5\qquad&&71\\ &E_7^{(1)}\qquad&&
p=1,5\qquad&&E_7\qquad&&A_1\times D_6\qquad&&125\\ &\qquad&&
p=7\qquad&&E_7\qquad&&A_7\qquad&&143\\ &E_8^{(1)}\qquad&&
p=1\qquad&&E_8\qquad&&A_1\times E_7\qquad&&239\\ &\qquad&&
p=7\qquad&&E_8\qquad&&D_8\qquad&&269
\endalignat
$$
\bigskip

If $k=2$ and $a_p=1$, we have

$$
\alignat5 &\text{type of $\widehat\g$}&&\qquad p &&\text{type of $\D_f$}\ &&\
\text{type of
$\g_0$}&&| \Cal W_{ab}^{\sigma} |\\ &A_{2n}^{(2)} \qquad&&
p=n\qquad&&C_n\qquad&&B_n\qquad&&2^{n+1}-1\\ &A_{2n-1}^{(2)}\qquad&& p=0,
1\qquad&&C_n\qquad&&C_n\qquad&&2^{n-1}\\ & &&
p=n\qquad&&C_n\qquad&&D_n\qquad&&2^{n+1}-1\\ &D_{n+1}^{(2)}\qquad&& 1\leq p\leq
n-1\qquad&&B_n\qquad&& B_p\times B_{n-p}\qquad&&4\binom{n}{p}-1\\ &\qquad&&
p=0,n\qquad&&B_n\qquad&&B_n\qquad&&2\\ &E_6^{(2)}\qquad&&
p=0\qquad&&F_4\qquad&&F_4\qquad&&4\\ &\qquad&& p=4\qquad&&F_4&&B_4\qquad&&23\\
\endalignat
$$

\bigskip

\heading \S6 The hermitian symmetric case.\endheading
\medskip

In this section we assume that $\g_0$ is not semisimple. Since
$k\, \sum\limits_{i=0}^n a_is_i=2$ this happens if and only if
$k=1$ and there are two indices $p,q$ such that
$a_p=a_q=s_p=s_q=1$ and $s_i=0$ for $i\ne p,q$.  By Theorem A (3)
we can and do choose $p=0$.
\par

Set $$D_\s'=D_\s\cap \{(x,\a_q)<0\}, \qquad
D_\s''=D_\s\cap \{(x,\a_q)\geq 0\}.$$
Clearly, $\Vol(D_\sigma)\big/\Vol(C_1)=\Vol(D'_\s)\big/\Vol(C_1)+
\Vol(D''_\s)\big/\Vol(C_1)$. We first compute
$\Vol(D'_\s)\big/\Vol(C_1)$.
\medskip

Denote by $P^\vee_f$, $Q^\vee_f$,  the coweight and the
coroot lattices of $\D_f$, and by $\ell_f$ its connection 
index, $\ell_f=[P_f:Q_f]$. 
\par
As in Section 5, let $\ovee_1,\dots,  \ovee_n$ be the fundamental
coweights of $\D_f$. Moreover let $\h^*_\s$ be the real span of
$\Da_0$. Then we have an orthogonal decomposition
$$
\h^*_0=\h_\s^*\oplus\R\omega_q^\vee;
$$
we denote by $\pi_\s$ to the corresponding projection onto $\h_\s^*$.
\par
It is clear that 
$\{\pi_\s(\ovee_i)\mid i\ne q\}$ is the dual basis of
$\Pia_0$ in $\h^*_\s$.
We denote by $P^\vee_\s$, $Q^\vee_\s$, $\ell_\s$  the coweight
lattice, the coroot lattice and the connection index of $\Da_0$:
$P_\s^\vee=\sum\limits_{i\ne 0,q}\Bbb Z\,\pi_\s(\ovee_i)$,
$Q^\vee_\s=\sum\limits_{i\ne 0, q}\Bbb Z\, \a^\vee_i$, and
$\ell_\s=[P_\s^\vee: Q_\s^\vee]$.

\proclaim{Lemma~6.1}
$$
{\Vol(D'_\s)\over\Vol(C_1)}={\ell_\s| W_f |\over \ell_f| W_\s |}.
$$
\endproclaim

\demo{Proof}
Set 
$$
\widehat I_\s=\{\sum_{i\ne q}x_i\omega^\vee_i\mid 0\leq x_i\leq 1\},
\quad I_q=\{x \ovee_q\mid 0\leq x\leq 1\}, \quad I=\widehat I_\s+I_q,
$$
and
$$
I'=\{\sum_{i\ne0}x_i\a^\vee_i\mid 0\leq x_i\leq 1\}.
$$
Then $\omega_0+I$ and $\omega_0+I'$ are fundamental domains for 
the action of $T(P_f^\vee)$ and $T(Q_f^\vee)$ on $\huno$,  
hence by \cite{1, VI.4, Lemma 1} we have that
$$
{\Vol(I)\over\Vol(C_1)}={\Vol(I)\over\Vol(I')}{\Vol(I')\over\Vol(C_1)}=
{1 \over \ell_f}| W_f |.
$$
\par
$\Wa_\s$ acts faithfully on $\omega_0+\h^*_\s$. Set 
$$I_\s=\pi_\s(\widehat I_\s); \quad
A_\si=(\bigcap\limits_{\a\in\Pia_\s}\ H_\a^+)\cap (\omega_0+\h^*_\s);
\quad 
\widehat A_\si=(\bigcap\limits_{\a\in\Pia_\s} H^+_\a)\cap
(\omega_0+ \sum\limits_{i\ne q}\real \ovee_i).$$
We notice that $A_\si=\pi_\s(\widehat A_\si),$ and that 
$D'_\s=\widehat  A_\si-I_q$.
\par
Arguing as in the proof of Lemma 5.4 we obtain that $A_\si$ is a
fundamental domain for the action of $\Wa_\si$ on
$\omega_0+\h^*_\s$. Moreover, $\omega_0+I_\s$ is a fundamental
domain for the action of $T(P_\s^\vee)$ onto $\omega_0+\h^*_\s$,
hence we also have that
$$
{\Vol_{n-1}(I_\s)\over\Vol_{n-1}(A_\s)}
={1\over \ell_\s}| W_\s |.
$$
Now we observe that  $\Vol(I)=\Vol(I_\s+I_q)=\Vol_{n-1}(I_\s)$,
and similarly $\Vol(D'_\s)=\Vol(A_\s-I_q)=\Vol_{n-1}(A_\s)$. It
follows that
$$
{\Vol(D'_\s)\over\Vol(C_1)}={\Vol(D'_\s)\over \Vol(I)}
{\Vol(I)\over \Vol(C_1)}=
{\Vol_{n-1}(A_\s)\over 
\Vol_{n-1}(I_\s)}{\Vol(I)\over\Vol(C_1)}=
{\ell_\s| W_f |\over \ell_f| W_\s |}.
$$
\endemo
\bigskip

We now compute $\Vol(D_\s'')\big/\Vol(C_1)$.
\proclaim{Lemma~6.2}
$${\Vol(D_\s'')\over \Vol(C_1)}={|W_f|\over |W_{\s}|}.$$
\endproclaim

\demo{Proof}
For $w\in W_f$ we denote by $\des(w)$ the descent set of $w$,
i.e. $\des(w)=N(w)\cap \Pi_f$.  For $I\subset \Pi_f$, set
$X^I=\{w\in W_f\mid \des(w)\cap I=\emptyset\}$. Note that
$W_\s=\langle s_\a\mid \a\in I\rangle$ with
$I=\Pi_f\setminus\{\a_q\}$. It is well known that $X^I$ is the
set of minimal length representatives of the right cosets
$W_\s\backslash W_f$. In particular $|X^I|=[W_f:W_\s]$. We shall
prove that $D_\s''-o_q=X_q C_1$, where $X_q=\{w\in W_f \mid
\des(w)\subseteq\{\a_q\}\}$, and $o_q={\ovee_q\over a_q}$; by the
above discussion this implies the claim.
\par
It is well known (see \cite{2, Lemma 1.2}) that $W_f C_1=
\{x\in\huno\mid -1\leq (\b,x)\leq 1\  \text {for all} \ \b\in
\D_f^+\}$. Moreover $\des(w) \subseteq \{\a_q\}$ if and only if
$(\a_i, w(C_1)) \geq 0$ for $i\ne q$. Hence $X_q
C_1=\{x\in\huno\mid -1\leq (\b,x) \leq 1\  \text {for all} \
\b\in \D_f^+, \ (\a_i, x)\geq 0\ \text{for} \ i\ne q\}$. 
\par
We first prove that if $x\in D_\s''$, then $x-o_q\in X_q C_1$. We
notice that $o_q=\ovee_q$, that $(\b,o_q)\leq 1$ for all $\b\in
\D_f^+$, and $(\b, o_q)=0$ if and only if $\b\in \Da_0$.
Moreover, $x$ is a dominant element, hence for positive $\b,
\a\in \D_f$ such that $\b\leq \a$ (in the standard partial order
on roots) we have $0\leq(\b,x)\leq(\a,x)$. Now assume that $\b \in
\D_f^+$. If $\b\in \Da_0^+$, then $\b\in \D_\Si$ for some
$\Si|\Pi_\s$, hence $0\leq (x,\b)=(x-o_q, \b) \leq (x,
\th_\Si)\leq 1$. If $\b\not\in \Da_0^+$, then $-1\leq
(x,\b)-1=(x-o_q, \b)\leq (x, \th)-1\leq 2-1=1$. 
\par
Next we prove the reverse inclusion. We consider $y\in X_q C_1$
and prove that $y+o_q\in D_\s''$. We have $(y+o_q, \a_q)=(y,
\a_q)+1\geq 0$, and if $i\ne q$ we have $(y+o_q, \a_i)=(y,
\a_i)\geq 0$. Moreover, $(y+o_q, \th)=(y, \th)+1\leq 2$, and
$(y+o_q, \th_\Si)=(y, \th_\Si)\leq 1$ for all $\Si|\Pi_\s$. This
concludes the proof.
\endemo

Combining the two lemmas we find

\proclaim{Theorem 6.3}If $\g_0$ is not semisimple then
$$ |\Cal  W_{ab}^{\sigma} |  =\frac{| W_f |}
{| W_\s |}\left(1+\frac{\ell_\s}{\ell_f}\right).
$$
\endproclaim
\bigskip

\ni 
We summarize the explicit results in the following table:

$$
\alignat5 &\text{type of $\widehat\g$}&&\qquad q &&\text{type of $\D_f$}\ &&\
\text{type of
$[\g_0,\g_0]$}&&| \Cal W_{ab}^{\sigma} |\\
&A_{n}^{(1)} \qquad&&
1\le q\le n\qquad&&A_n\qquad&&A_{q-1}\times A_{n-q}\qquad&&\binom{n+1}{q}+q\binom{n}{q}\\
&B_{n}^{(1)} \qquad&&
q=1\qquad&&B_n\qquad&&B_{n-1} \qquad&&4n\\
&C_{n}^{(1)} \qquad&&
q=n\qquad&&C_n\qquad&&A_{n-1} \qquad&&2^{n-1}(n+2)\\
&D_{n}^{(1)} \qquad&&
q=1\qquad&&D_n\qquad&&D_{n-1} \qquad&&4n\\
&\qquad&&
q=n-1,n\qquad&&D_n\qquad&&A_{n-1} \qquad&&2^{n-3}(n+4)\\
&E_{6}^{(1)} \qquad&&
q=1,6\qquad&&E_6\qquad&&D_{5} \qquad&&63\\
&E_{7}^{(1)} \qquad&&
q=7\qquad&&E_7\qquad&&E_{6} \qquad&&140
\endalignat
$$

\enddocument
\bye